# STATIONARY DISTRIBUTION FOR DIOECIOUS BRANCHING PARTICLE SYSTEMS WITH RAPID STIRRING

By Feng Yu[1]

*University of Oxford*

We study dioecious (i.e., two-sex) branching particle system models, where there are two types of particles, modeling the male and female populations, and where birth of new particles requires the presence of both male and female particles. We show that stationary distributions of various dioecious branching particle models are nontrivial under certain conditions, for example, when there is sufficiently fast stirring.

**1. The particle models.** We consider a type of particle system that can be used to model sexual reproduction of a certain species. This work was inspired in part by Dawson and Perkins [1], which studied the following system of stochastic partial differential equations:

$$\frac{\partial u}{\partial t}(t,x) = \frac{1}{2}\Delta u(t,x) + (\gamma u(t,x)v(t,x))^{1/2}\dot{W}_1(t,x),$$
$$\frac{\partial v}{\partial t}(t,x) = \frac{1}{2}\Delta v(t,x) + (\gamma u(t,x)v(t,x))^{1/2}\dot{W}_2(t,x),$$

where $\Delta = \sum_i \partial^2/\partial x_i^2$ is the Laplacian, $\gamma > 0$ and $\dot{W}_i(t,x)$ ($i=1,2$) are independent space-time white noises on $\mathbb{R}_+ \times \mathbb{R}$. One can associate $u(x,t)$ and $v(x,t)$ with the male and female populations of particles, respectively, at spatial location $x$ and time $t$. Loosely speaking, the above SPDE system says that individual male or female particles move around according to Brownian motion, but branching is only possible when both male and female particles are present at the same spatial location. Notice that at spatial locations

Received October 2005; revised March 2007.

[1]Supported by the Department of Mathematics and a Killam Fellowship while at the University of British Columbia, Vancouver, Canada, and EPSRC Grant GR/T19537 while at the University of Oxford, UK.

where the female population is 0, the branching rate for the male population is also 0, therefore the male population does not die and the only effect on the male population at those spatial locations is the diffusive effect of the heat kernel. This behavior is not very realistic, since one would expect a natural death rate for male particles, even without the presence of any female particles. In this work, we study models involving a finite number of male and female particles with more realistic behavior.

The model we study involves two types of particles, male and female, residing on the grid $S = \mathbb{Z}^d$. Each *site* $x \in S$ contains two *nests*, one for the male particle and the other for the female particle. Each nest can be inhabited by at most one particle, either male or female. Let $E = \{0,1\}$ and $F = E \times E$ be the set of possible states at each site in $S$. For $x \in S$, we write

$$\xi(x) = (\xi^1(x), \xi^2(x)),$$

where $\xi^1(x)$ denotes the number (0 or 1) of male particles at site $x$ and $\xi^2(x)$ denotes the number of female particles at site $x$. We define the interaction neighborhood

$$\mathcal{N} = \{0, y_1, \ldots, y_N\}$$

and the neighborhood of $x$

$$\mathcal{N}_x = x + \mathcal{N}.$$

For example, $\mathcal{N} = \{0, -1, 1\}$ if the interaction is nearest-neighbor on $\mathbb{Z}$. Let $c_i(x, m, \xi)$ denote the rate at which nest $m$ ($m = 1, 2$) of site $x$ flips to state $i$ ($i = 0, 1$) and assume that $c_i(x, m, \xi)$ depends only on the neighborhood $\mathcal{N}_x$, that is,

$$c_i(x, m, \xi) = h_{i,m}(\xi(x), \xi(x+y_1), \ldots, \xi(x+y_N))$$

for some function $h_{i,m} : F^{N+1} \to \mathbb{R}^+$. The death rate $c_0$ is always constant,

$$(1) \qquad c_0(x, m, \xi) = \begin{cases} \delta, & \text{if } \xi^m(x) = 1, \\ 0, & \text{otherwise,} \end{cases}$$

while the birth rate $c_1(x, m, \xi)$ is positive only if both male and female particles can be found in $\mathcal{N}_x$. In this work, we take $\delta = 1$. Note that this simply means that the unit of time we take is the average lifetime of an individual. For example, the dioecious branching particle model which we will consider in Section 1.1 has

$$(2) \qquad c_1(x, m, \xi) = \begin{cases} \lambda n_1(x, \xi) n_2(x, \xi), & \text{if } \xi^m(x) = 0, \\ 0, & \text{otherwise,} \end{cases}$$

where

$$n_{m'}(x, \xi) = |\{z \in \mathcal{N}_x : \xi^{m'}(x+z) = 1\}|,$$

that is, at rate $\lambda$, each pair of male and female particles in $\mathcal{N}_x$ give birth to a particle at nest $m$ of site $x$ if that nest is not already occupied. A more stringent condition, as in the particle model with rapid stirring which we will consider in Section 1.2, is to require both parent particles to reside at the same site, that is,

$$(3) \qquad c_1(x, m, \xi) = \begin{cases} \lambda n_{1+2}(x, \xi), & \text{if } \xi^m(x) = 0, \\ 0, & \text{otherwise,} \end{cases}$$

where

$$n_{1+2}(x, \xi) = |\{z \in \mathcal{N}_x : \xi^1(x+z) = 1 \text{ and } \xi^2(x+z) = 1\}|.$$

This more stringent condition should not alter the behavior of the particle system if one allows a larger $\lambda$ than in (2), but it does help to simplify the analysis somewhat.

1.1. *Dioecious branching particle model.* We first describe the model with birth and death rates as in (1) and (2), for which we will establish the existence of nontrivial stationary distribution(s) and, consequently, a phase transition (Section 2). The birth and death rates at site $x$ depend only on $\xi$ in a neighborhood of $x$, therefore birth (and death) rates at $x$ and $y$ where $\mathcal{N}_x \cap \mathcal{N}_y = \varnothing$ are independent. First, we restate the model:

1. *Birth.* For each nest $(x, m)$ and each pair $(z_1, z_2) \in \mathcal{N}_x \times \mathcal{N}_x$ such that $\xi^1(z_1) = 1$ and $\xi^2(z_2) = 1$, where $z_1$ and $z_2$ need not be distinct, with rate $\lambda$, a child of $(z_1, z_2)$ is born into nest $m$ of site $x$ if $(x, m)$ is not already occupied.
2. *Death.* Each particle dies at rate 1.

We can think of this particle system as a generalized spin system, generalized in the sense that the phase space at each site is $\{0,1\}^2$ rather than $\{0,1\}$. One can refer to Chapter 3 of Liggett [5] for a detailed introduction to classic spin systems. We observe that the all-0 state [i.e., $\xi^1(x) = \xi^2(x) = 0$ for all $x$] is an absorbing state, therefore the probability measure that concentrates only on the all-0 state is a trivial stationary distribution. We say that a stationary distribution is *nontrivial* if it does not concentrate only on the all-0 state. A major goal of this work is to establish the existence of nontrivial stationary distributions for various particle systems.

The interacting particle system involving the birth and death mechanisms described above can be constructed using a countable number of Poisson processes, similar to the construction found in Chapter 2 of Durrett [2]. Define

$$c^* = \sup_{\xi, m} \sum_i c_i(x, m, \xi).$$

We assume $c^* < \infty$. Let $\{T_n^{x,i,m} : n \geq 1\}$ be the arrival times of independent rate $c^*$ Poisson processes and $\{U_n^{x,i,m} : n \geq 1\}$ be independent uniform random variables on $[0,1]$. At time $t = T_n^{x,i,m}$, nest $(x,m)$ flips to state $i$ if $U_n^{x,i,m} \leq c_i(x,m,\xi_{t-})/c^*$ and stays unchanged otherwise.

Alternatively, one can explicitly write down the generators $\mathcal{G}^1$ and $\mathcal{G}^2$ associated with the particle system with death rates (1) and birth rates (2) and (3), respectively, as follows:

$$\mathcal{G}^1 f(\xi) = \sum_{(x,m) \in S \times \{1,2\}} \left[ \xi^m(x)(f(\xi - \delta_{x,m}) - f(\xi)) \right.$$
$$+ \sum_{y,z \in \mathcal{N}_x} \lambda \xi^1(y) \xi^2(z)(1 - \xi^m(x))$$
$$\left. \times (f(\xi + \delta_{x,m}) - f(\xi)) \right] \quad (4)$$

and

$$\mathcal{G}^2 f(\xi) = \sum_{(x,m) \in S \times \{1,2\}} \left[ \xi^m(x)(f(\xi - \delta_{x,m}) - f(\xi)) \right.$$
$$+ \sum_{y \in \mathcal{N}_x} \lambda \xi^1(y) \xi^2(y)(1 - \xi^m(x))$$
$$\left. \times (f(\xi + \delta_{x,m}) - f(\xi)) \right], \quad (5)$$

where $f$ is a function on $F^S$ (endowed with the product topology) with compact support (in this case, compact support implies that $f$ depends on only finitely many sites) and $\delta_{x,m}$ is a function on $S \times \{1,2\}$ that is one at $(x,m)$ and zero elsewhere, and apply Theorem B3 in Liggett [6] to see that the closure of $\mathcal{G}$ [in $C(F^S) \times C(F^S)$, with respect to the sup-norm topology] generates a Feller Markov process. Summarizing results from the two preceding paragraphs, we have the following theorem.

THEOREM 1.1. *There exists a unique Feller Markov process $\xi_t$ constructed as before with generator* (4) *or* (5).

The particle system $\xi$ with generator (4) or (5) is *attractive* in the sense that $\xi$ is monotonic in initial conditions. One can check that if $\xi_0(x) \leq \bar{\xi}_0(x)$ for all $x \in S$, where $\leq$ denotes the partial order $(0,0) \leq (0,1), (1,0) \leq (1,1)$, then $\xi_t(x) \leq \bar{\xi}_t(x)$ for all $x$ and $t$. This is true since every birth or death event preserves $\leq$. For example, if $\xi_{t-}(x) = (0,0)$, $\bar{\xi}_{t-}(x) = (0,1)$ and at time $t$ there is a male birth event at site $x$, then $\xi_t(x) = (1,0)$ and $\bar{\xi}_t(x) = (1,1)$, so the inequality $\xi_t(x) \leq \bar{\xi}_t(x)$ has been maintained. Similarly, one can check that the particle system $\xi$ is increasing in the birth rate $\lambda$ by coupling the random variables $T_n^{x,i,m}$ and $U_n^{x,i,m}$ involved in the constructions in the obvious way. Because of this monotonicity, along with the existence of nontrivial stationary distributions for sufficiently large $\lambda$ and extinction for sufficiently small $\lambda$ which we will establish later in this work, we may conclude that there is a phase transition in the behavior of the particle system $\xi$.

1.2. *Description of the particle model with rapid stirring.* If we add rapid stirring to the particle system, that is, we scale the integer grid $\mathbb{Z}^d$ by $\varepsilon$ and stir neighboring particles at rate $\varepsilon^{-2}$ in addition to performing the birth and death mechanisms, then the particle system converges to the solution of a reaction-diffusion PDE as $\varepsilon \to 0$ (see Theorems 8.1 and 8.2 in Durrett [2] and the beginning of Section 1.3 of this work). This PDE represents the mean-field behavior of the particle system and is usually easier to analyze than the particle system itself. As promised earlier, we will establish in Section 2 that there is a phase transition for the dioecious branching particle model (i.e., without rapid stirring), but obtaining any reasonable estimates on exactly where this transition occurs seems to be difficult. One advantage of adding rapid stirring mechanisms is that one can get a reasonably good idea of where the phase transition occurs in the rapidly stirred particle model by analyzing the limiting PDE, or simulating this PDE on a computer.

Moreover, this convergence establishes a connection between the particle model and PDE systems, which is of independent interest. Since many PDE's arise out of natural systems, this connection justifies the study of the PDE. The underlying stochastic system can also yield information about the PDE; for example, as we will see in Section 1.3, the monotonicity of the particle system will lead to the monotonicity of the PDE. Information about the PDE will similarly yield information about the particle model. In Section 4, we will establish a condition on the PDE that implies the existence of nontrivial stationary distributions for the particle system with sufficiently small $\varepsilon$.

For the particle models with rapid stirring, we work with $S = \varepsilon \mathbb{Z}^d$ and denote the corresponding process by $\xi^\varepsilon$. We also assume the birth rates in (3) and death rate $\delta = 1$, while the neighborhood $\mathcal{N}$ is nearest-neighbor: $\mathcal{N} = \{y : \|y\| = 0 \text{ or } \varepsilon\}$. Here, we use the $L^1$-norm: $\|y\| = \sum_{k=1}^d |y_k|$. In addition to the transitions in the dioecious branching model, we introduce spatial movement of particles between neighboring sites called *rapid stirring*. We consider two rapid stirring mechanisms in this work, one called *lily-pad* stirring and the other called *individual* stirring:

- *Lily-pad stirring.* For each $x, y \in \varepsilon \mathbb{Z}^d$ with $\|x - y\|_1 = \varepsilon$, $\xi^\varepsilon(x) = (\xi^{\varepsilon,1}(x), \xi^{\varepsilon,2}(x))$ and $\xi^\varepsilon(y) = (\xi^{\varepsilon,1}(y), \xi^{\varepsilon,2}(y))$ are exchanged at rate $\varepsilon^{-2}$.

- *Individual stirring.* For each $i \in \{1,2\}$ and $x, y \in \varepsilon \mathbb{Z}^d$ with $\|x - y\|_1 = \varepsilon$, $\xi^{\varepsilon,i}(x)$ and $\xi^{\varepsilon,i}(y)$ are exchanged at rate $\varepsilon^{-2}$.

Just as in the particle model without rapid stirring described in Section 1.1, one can construct the particle model with either lily-pad stirring or individual stirring using a countable number of Poisson processes. Alternatively, one can write down the generator explicitly and again apply Theorem B3 in Liggett [6] to establish:

THEOREM 1.2. *Let $S = \varepsilon \mathbb{Z}^d$ and $\mathcal{N} = \{y : \|y\| = 0 \text{ or } \varepsilon\}$. There exists a unique Feller process $\xi_t$ with generator $\mathcal{G}^L$ for the particle model with lily-pad stirring or generator $\mathcal{G}^I$ for the particle model with individual stirring,*

(6) $\quad \mathcal{G}^L f(\xi) = \mathcal{G}^2 f(\xi) + \sum_{x,y \in S, x \in \mathcal{N}_y} \varepsilon^{-2}(f(\xi^{x \leftrightarrow y}) - f(\xi))$,

(7) $\quad \mathcal{G}^I f(\xi) = \mathcal{G}^2 f(\xi) + \sum_{m \in \{1,2\}, x,y \in S, x \in \mathcal{N}_y} \varepsilon^{-2}(f(\xi^{(x,m) \leftrightarrow (y,m)}) - f(\xi))$,

*where*

$$\xi^{x \leftrightarrow y}(z, m') = \begin{cases} \xi(z, m'), & \text{if } z \neq x, y, \\ \xi(x, m'), & \text{if } z = y, \\ \xi(y, m'), & \text{if } z = x \end{cases}$$

*and*

$$\xi^{(x,m) \leftrightarrow (y,m)}(z, m') = \begin{cases} \xi(z, m'), & \text{if } (z, m') \neq (x, m), (y, m), \\ \xi(x, m), & \text{if } (z, m') = (y, m), \\ \xi(y, m), & \text{if } (z, m') = (x, m). \end{cases}$$

For lily-pad stirring, instead of thinking of a site that consists of two nests, as in the dioecious branching model, we can view each site as having for states in

$$F = \{0,1\}^2 = \{(0,0), (0,1), (1,0), (1,1)\}.$$

We restate the dynamics of the particle model in terms of these four states. At any site $x \in \varepsilon \mathbb{Z}^d$, only the following transitions are possible: $(0,0) \leftrightarrow (0,1)$, $(0,1) \leftrightarrow (1,1)$, $(0,0) \leftrightarrow (1,0)$ and $(1,0) \leftrightarrow (1,1)$, that is, only one particle is born or dies at a particular time. The rates of these transitions are as follows:

$c_{(0,0)}(x, \xi^\varepsilon) = 1,$ if $\xi^\varepsilon(x) = (0,1)$ or $\xi^\varepsilon(x) = (1,0)$,

$c_{(0,1)}(x, \xi^\varepsilon) = c_{(1,0)}(x, \xi^\varepsilon) = 1,$ if $\xi^\varepsilon(x) = (1,1)$,

$c_{(0,1)}(x, \xi^\varepsilon) = c_{(1,0)}(x, \xi^\varepsilon) = \lambda n_{1+2}(x, \xi^\varepsilon),$

if $\xi^\varepsilon(x) = (0,0)$,

$c_{(1,1)}(x, \xi^\varepsilon) = \lambda n_{1+2}(x, \xi^\varepsilon),$ if $\xi^\varepsilon(x) = (0,1)$ or $\xi^\varepsilon(x) = (1,0)$.

The difference between these two stirring mechanisms is that lily-pad stirring forces male and female particles at a site to move together, but individual stirring allows independent movement of male and female particles. Every exchange of particles, in both lily-pad stirring and individual stirring, is monotonicity preserving, thus neither stirring mechanism disrupts the monotonicity property of the particle system.

1.3. *Convergence to a PDE for lily-pad stirring.* Consider the particle system with lily-pad stirring and its generator given by (6). For $i \in F$, if we define

$$u_i^\varepsilon(t,x) = P(\xi_t^\varepsilon(x) = i),$$

then Theorem 8.1 in Durrett [2] (or its generalization, Theorem 3.1 of this work) shows that if $g_i(x)$ is $C^1$ and $u_i^\varepsilon(0,x) = g_i(x)$ for all $i$, then

$$(8) \qquad u_i(t,x) = \lim_{\varepsilon \to 0} u_i^\varepsilon(t,x)$$

exists and satisfies the following system of PDEs:

$$(9) \quad \begin{aligned} \frac{\partial u_{(0,0)}}{\partial t} &= \Delta u_{(0,0)} + u_{(0,1)} + u_{(1,0)} - 4\lambda d u_{(0,0)} u_{(1,1)}, \\ \frac{\partial u_{(0,1)}}{\partial t} &= \Delta u_{(0,1)} + u_{(1,1)} - u_{(0,1)} + 2\lambda d (u_{(0,0)} - u_{(0,1)}) u_{(1,1)}, \\ \frac{\partial u_{(1,0)}}{\partial t} &= \Delta u_{(1,0)} + u_{(1,1)} - u_{(1,0)} + 2\lambda d (u_{(0,0)} - u_{(1,0)}) u_{(1,1)}, \\ \frac{\partial u_{(1,1)}}{\partial t} &= \Delta u_{(1,1)} - 2u_{(1,1)} + 2\lambda d (u_{(0,1)} + u_{(1,0)}) u_{(1,1)}. \end{aligned}$$

In order for the limit in (8) to make sense, we extend the definition of $u_i^\varepsilon(t,x)$ to all $x \in \mathbb{R}^d$ by requiring $u_i^\varepsilon(t,x) = u_i^\varepsilon(t, \varepsilon \lfloor x/\varepsilon \rfloor)$ where, $\lfloor q \rfloor$ denotes the integer part of $q$. Obviously, $u_i$ must lie in $[0,1]$ for all $i$ and $t$ since it is a limit of probabilities. The reaction terms of the equation involving $\partial u_{(0,0)}/\partial t$ say roughly that sites enter state $(0,0)$ at rate $u_{(0,1)} + u_{(1,0)}$, that is, when the only male or female particle at the site dies; sites leave state $(0,0)$ at rate $4\lambda d u_{(0,0)} u_{(1,1)}$, that is, when a particular site has no particles and a particle (male or female) is born by a pair of male and female particles at one of the $2d$ neighboring sites.

We also observe that, strictly speaking, one should require the initial condition of (9) to be $C^2$, for otherwise, the PDE system may not make sense at $t=0$. This issue can be remedied by considering the corresponding integral equation, as in equation (15.1.2) of Taylor [7]. Furthermore, Proposition 15.1.2 of Taylor [7] shows that solutions to (9) are $C^\infty$ at any $t > 0$ as long as the initial condition is $C^1$. Thus we only require our initial conditions to be $C^1$ from this point on.

We transform the parameter space of the 3-dimensional system in (9) to obtain a monotone 2-dimensional system, which will be easier to analyze. First, define $v_0 = u_{(0,0)}$, $v_1 = u_{(0,1)} + u_{(1,0)}$ and $v_2 = u_{(1,1)}$, $(v_0, v_1, v_2)$ then satisfies

$$\frac{\partial v_0}{\partial t} = \Delta v_0 + v_1 - 4\lambda d v_0 v_2,$$

$$\frac{\partial v_1}{\partial t} = \Delta v_1 + 2v_2 - v_1 + 2\lambda d(2v_0 - v_1)v_2,$$

$$\frac{\partial v_2}{\partial t} = \Delta v_2 - 2v_2 + 2\lambda d v_1 v_2,$$

where $v_0 + v_1 + v_2 = 1$. The above system can be written as the limiting PDE under rapid stirring of another particle system $\zeta^\varepsilon$, still on $S = \varepsilon \mathbb{Z}^d$, with state space $F = \{0, 1, 2\}$ and transitions $0 \leftrightarrow 1$ and $1 \leftrightarrow 2$ at rates

$$\begin{aligned} c_0(x, \zeta^\varepsilon) &= 1, & \text{if } \zeta^\varepsilon(x) = 1, \\ c_1(x, \zeta^\varepsilon) &= 2, & \text{if } \zeta^\varepsilon(x) = 2, \\ c_1(x, \zeta^\varepsilon) &= 2\lambda n_2(x, \zeta^\varepsilon), & \text{if } \zeta^\varepsilon(x) = 0, \\ c_2(x, \zeta^\varepsilon) &= \lambda n_2(x, \zeta^\varepsilon), & \text{if } \zeta^\varepsilon(x) = 1, \end{aligned}$$

where

$$n_2(x, \xi^\varepsilon) = |\{z \in \mathcal{N} : \zeta^\varepsilon(x + z) = 2\}|.$$

One can interpret $\zeta^\varepsilon(x)$ as the total number of $\xi^\varepsilon$-particles at $x$ and check that monotonicity in initial condition holds for $\zeta^\varepsilon$. Let $(v_1^\varepsilon(0, x), v_2^\varepsilon(0, x)) = (g_1(x), g_2(x))$ and $(\bar{v}_1^\varepsilon(0, x), \bar{v}_2^\varepsilon(0, x)) = (\bar{g}_1(x), \bar{g}_2(x))$ be two sets of initial distributions such that $g_2 \leq \bar{g}_2$ and $g_1 + g_2 \leq \bar{g}_1 + \bar{g}_2$, with $g_1, g_2, g_1 + g_2, \bar{g}_1, \bar{g}_2, \bar{g}_1 + \bar{g}_2$ all lying in $[0, 1]$. It is then possible to set up two initial conditions $\zeta_0^\varepsilon$ and $\bar{\zeta}_0^\varepsilon$ such that $P(\zeta_0^\varepsilon(x) = i) = v_i^\varepsilon(0, x)$ and $P(\bar{\zeta}_0^\varepsilon(x) = i) = \bar{v}_i^\varepsilon(0, x)$, $i = 1, 2$, and $\zeta_0^\varepsilon(x) \leq \bar{\zeta}_0^\varepsilon(x)$ holds for all $x$ and $\omega$. The monotonicity property of $\zeta^\varepsilon$ implies that $\zeta_t^\varepsilon(x) \leq \bar{\zeta}_t^\varepsilon(x)$ for all $t$ and $x$ and therefore

$$P(\zeta_t^\varepsilon(x) \geq 1) \leq P(\bar{\zeta}_t^\varepsilon(x) \geq 1) \quad \text{and} \quad P(\zeta_t^\varepsilon(x) \geq 2) \leq P(\bar{\zeta}_t^\varepsilon(x) \geq 2),$$

that is, for all $t$ and $x$,

$$v_2^\varepsilon(t, x) \leq \bar{v}_2^\varepsilon(t, x),$$

$$v_1^\varepsilon(t, x) + v_2^\varepsilon(t, x) \leq \bar{v}_1^\varepsilon(t, x) + \bar{v}_2^\varepsilon(t, x).$$

We now transform the parameter space a second time by defining $(u, v) = (v_1 + v_2, v_2)$ and writing $c = 2\lambda d$. $(u_t, v_t)$ is then monotone in initial conditions since $v_i^\varepsilon(t, x) \to v_i(t, x)$. We observe that $u$ can be interpreted as the density of occupied sites (where either one or both nests are occupied) and $v$ as the density of doubly occupied sites (where both nests are occupied). Straightforward calculation shows that $(u,v)$ satisfies the following system:

$$\begin{aligned}\frac{\partial u}{\partial t} &= \Delta u + (2c(1-u)+1)v - u, \\ \frac{\partial v}{\partial t} &= \Delta v + (c(u-v)-2)v.\end{aligned} \tag{10}$$

We summarize this paragraph in the following lemma.

LEMMA 1.3. *The PDE system* (10) *is monotone in initial conditions that lie in*

$$\mathcal{R} = \{(u,v): 0 \leq v \leq u \leq 1\},$$

*that is, if there are two initial conditions* $(u_0, v_0) \in \mathcal{R}$ *and* $(\bar{u}_0, \bar{v}_0) \in \mathcal{R}$, *with* $u_0 \leq \bar{u}_0$ *and* $v_0 \leq \bar{v}_0$ *everywhere, then* $u_t \leq \bar{u}_t$ *and* $v_t \leq \bar{v}_t$ *everywhere, for all* $t$. *Furthermore, both* $(u_t, v_t)$ *and* $(\bar{u}_t, \bar{v}_t)$ *lie in* $\mathcal{R}$ *for all* $t$.

In Section 4, we will analyze (10) to establish the following result.

THEOREM 1.4. *If* $\lambda$ *is sufficiently large and* $\varepsilon$ *is sufficiently small, then there exists a nontrivial translation invariant stationary distribution for the dioecious branching particle model with lily-pad stirring with generator* (6).

1.4. *Convergence to a PDE for individual stirring.* We consider the particle system with individual stirring and with its generator given by (7). Unlike lily-pad stirring, Theorem 8.1 in Durrett [2] cannot be directly applied to obtain convergence to a PDE system for individual stirring. We can, however, follow the ideas used in the proof of that theorem to establish a corresponding result, Theorem 3.1. For the process $\xi^\varepsilon$ with generator (7) and $i \in E$, define

$$u^\varepsilon_{i,m}(t,x) = P(\xi^\varepsilon_t(x,m) = i).$$

Theorem 3.1 then implies that if $g_{i,m}: \mathbb{R} \to [0,1]$ is $C^1$ and $u^\varepsilon_{i,m}(0,x) = g_{i,m}(x)$, then $u_{i,m}(t,x) = \lim_{\varepsilon \to 0} u^\varepsilon_{i,m}(t,x)$ exists and satisfies the following system of PDEs:

$$\begin{aligned}\frac{\partial u_{1,1}}{\partial t} &= \Delta u_{1,1} - u_{1,1} + 2\lambda d(1-u_{1,1})u_{1,1}u_{1,2}, \\ \frac{\partial u_{1,2}}{\partial t} &= \Delta u_{1,2} - u_{1,2} + 2\lambda d(1-u_{1,2})u_{1,1}u_{1,2},\end{aligned}$$

where $u_{0,1} + u_{1,1} = u_{0,2} + u_{1,2} = 1$. Notice that if we start with a symmetric initial condition, that is, $g_{i,1} = g_{i,2}$, then the solution to the above PDE is also symmetric. And if we define $u = u_{1,1} = u_{1,2}$, then we obtain the following PDE for $u$:

$$\tag{11} \frac{\partial u}{\partial t} = \Delta u + f(u), \qquad f(u) = -u + 2\lambda d(1-u)u^2.$$

This PDE has been analyzed in Durrett and Neuhauser [3] as their sexual reproduction model (Example 3 on page 291). In fact, it is not difficult to see that if $u_{1,1} = u_{1,2}$, then choosing the "father" from the male population is exactly the same as choosing the "father" from the female population, hence it is quite natural for this reduction to occur. Theorem 4 of Durrett and Neuhauser [3] states that if $2\lambda d > 4.5$ and $\varepsilon$ is sufficiently small, then their sexual reproduction model has nontrivial stationary distribution(s). Although this theorem does not directly apply to our particle system $\xi^\varepsilon$ with two types of particles because of the difference in stirring mechanisms, one can nevertheless work through the proof of Lemma 3.3 of Durrett and Neuhauser [3], while making obvious changes, to establish the following, similar, result.

- Let $0 < \rho_1 < \rho_0 < 1$ be the two nonzero roots of $f(u)$. Define $\beta = (\rho_0 - \rho_1)/10$ and $I_k = 2Lke_1 + [-L, L)^d$. If $\varepsilon$ is small, $L$ is large, and $\xi^\varepsilon(0)$ has density at least $\rho_1 + \beta$ of both male particles and female particles in $I_0$, then for sufficiently large $T$, with high probability, $\xi^\varepsilon(T)$ will have density at least $\rho_0 - \beta$ in $I_1$ and $I_{-1}$.

This result can then be fed into a comparison argument, comparing the particle system with oriented percolation, as on page 312 of Durrett and Neuhauser [3] or in the proof of Theorem 4.3 in Durrett [2], to establish the existence of nontrivial stationary distribution(s) for the particle system $\xi^\varepsilon$ under individual stirring with sufficiently small $\varepsilon$. We then have the following theorem.

THEOREM 1.5. *If $\lambda$ is sufficiently large and $\varepsilon$ is sufficiently small, then there exists a nontrivial translation invariant stationary distribution for the dioecious branching particle model with generator* (7).

1.5. *Discussion.* We can numerically solve PDEs using standard finite difference methods to obtain rough values of $\lambda_c$ for $d = 2$, where $\lambda_c$ is, as in Durrett and Neuhauser [3], the critical value of $\lambda$ necessary for longterm survival in the rapid stirring limit ($\varepsilon \to 0$). We recall that (10) and (11) are the rapid stirring limits of the particle systems with generators $\mathcal{G}^L$ and $\mathcal{G}^I$ in (6) and (7), respectively, both with birth-death mechanisms described by $\mathcal{G}^2$ that require parents to be at the same site. We can replace $\mathcal{G}^2$ with $\mathcal{G}^1$ (which only requires one parent to be in the neighborhood of the other parent) in (4) and obtain slightly different generators $\tilde{\mathcal{G}}^L f(\xi) = \mathcal{G}^1 f(\xi) + \sum_{x,y \in S, x \in \mathcal{N}_y} \varepsilon^{-2}(f(\xi^{x \leftrightarrow y}) - f(\xi))$ and $\tilde{\mathcal{G}}^I f(\xi) = \mathcal{G}^1 f(\xi) + \sum_{m \in \{1,2\}, x,y \in S, x \in \mathcal{N}_y} \varepsilon^{-2}(f(\xi^{(x,m) \leftrightarrow (y,m)}) - f(\xi))$. The rapid stirring limit of $\tilde{\mathcal{G}}^I$ is exactly the PDE in (11) but with $2d$ replaced by $2d(2d+1)$ because there are $2d(2d+1)$ "potential parent-nest pairs" involved in $\mathcal{G}^1$. The rapid stirring limit of $\tilde{\mathcal{G}}^L$ satisfies the following PDE:

$$\frac{\partial u_{(0,0)}}{\partial t} = \Delta u_{(0,0)} - 2\lambda(2d)^2 u_{(0,0)}(u_{(0,1)} + u_{(1,1)})(u_{(1,0)} + u_{(1,1)}),$$

$$\frac{\partial u_{(0,1)}}{\partial t} = \Delta u_{(0,1)} - \lambda(2d)(2d+1)u_{(0,1)}(u_{(0,1)} + u_{(1,1)})(u_{(1,0)} + u_{(1,1)})$$
$$+ \lambda(2d)^2 u_0(u_{(0,1)} + u_{(1,1)})(u_{(1,0)} + u_{(1,1)}),$$

Table 1

|  | Individual stirring | Lily-pad stirring |
| --- | --- | --- |
| $\mathcal{G}^1$ | $\lambda_c = 0.225$ | $\lambda_c \in [0.271, 0.272]$ |
| $\mathcal{G}^2$ | $\lambda_c = 1.125$ | $\lambda_c \in [1.114, 1.115]$ |

(12) $\quad \frac{\partial u_{(1,0)}}{\partial t} = \Delta u_{(1,0)} - \lambda(2d)(2d+1)u_{(1,0)}(u_{(0,1)} + u_{(1,1)})(u_{(1,0)} + u_{(1,1)})$
$$+ \lambda(2d)^2 u_0(u_{(0,1)} + u_{(1,1)})(u_{(1,0)} + u_{(1,1)}),$$

$$\frac{\partial u_{(1,1)}}{\partial t} = \Delta u_{(1,1)} - \lambda(2d)(2d+1)(u_{(0,1)} + u_{(1,0)})$$
$$\times (u_{(0,1)} + u_{(1,1)})(u_{(1,0)} + u_{(1,1)}).$$

For individual stirring, the limiting PDE is (11) and Theorem 4 from Durrett and Neuhauser [3] shows that $\lambda_c = 4.5/2d = 1.125$ for $\mathcal{G}^I$ and $\lambda_c = 4.5/2d(2d+1) = 0.225$ for $\tilde{\mathcal{G}}^I$. For lily-pad stirring, we do not have a readily available theorem to tell us the exact value of $\lambda_c$. We obtained the range of values of $\lambda_c$ in the rightmost column of Table 1 by simulating (10) and (12) using progressively finer grids and stopping once changes in the estimates of $\lambda_c$ become smaller than $10^{-3}$.

Curiously, free movement of all individuals (individual stirring) seems to work better (to ensure survival at least) when mating occurs between all pairs of individuals in the neighborhood, as in $\mathcal{G}^1$, but restricting male and female individuals to move together (lily-pad stirring) seems to work better when mating occurs only between male and female individuals at the same site, as in $\mathcal{G}^2$. A possible explanation is that matching of movement of individuals and the mating strategy helps survival: when the mating strategy is to mate with any individual of the opposite sex in the neighborhood, free movement of all individuals helps everyone to find mating partners more easily; but when the mating strategy is to mate with only individuals at the same site, free movement of all individuals only helps to break up "marriages."

In the single-sex scenario where the birth and death mechanisms are the same as those of the contact process with rate of infection $\lambda$ and rate of recovery 1, we obtain the following PDE in the rapid stirring limit:

$$\frac{\partial u}{\partial t} = \Delta u + f(u), \qquad f(u) = -u + 2\lambda d(1-u)u,$$

which has $\lambda_c = 1/2d$ (i.e., $\lambda_c = 0.25$) by Theorem 2 from Durrett and Neuhauser [3]. This $\lambda_c$ is roughly the same as the $\lambda_c$ for birth-death mechanisms $\mathcal{G}^1$, although much smaller than the $\lambda_c$ for birth-death mechanisms $\mathcal{G}^2$ since the birth mechanism in $\mathcal{G}^2$ is much more restrictive. With birth-death mechanisms $\mathcal{G}^1$ and $d=2$, the density of individuals (male or female) at equilibrium when $\lambda$ is only slightly higher than the $\lambda_c$ is roughly 0.66 for individual stirring and 0.79 for lily-pad stirring. Thus, in the two-sex scenario $\mathcal{G}^1$, although each individual can only give birth if there is any individual of the opposite sex in the neighborhood, finding a mate in the neighborhood should not prove to be a problem and it is not terribly surprising that the $\lambda_c$ for $\mathcal{G}^1$ with individual stirring is smaller than the $\lambda_c$ for the single-sex scenario.

In the remainder of this paper, we will establish various results as promised in this section. In Section 2, we prove a few results on the dioecious branching particle system without rapid stirring. In Section 3, we present, without proof, a convergence theorem for particle system with individual stirring. Finally, in Section 4, we prove Theorem 1.4 by establishing a condition on (10), similar to condition (*) on page 180 of Durrett [2].

**2. Results on the dioecious branching particle model.** In this section, we assume the model with generator (4) where we take $S = \mathbb{Z}^d$, that is, the particle system with birth and death mechanisms, but no stirring. We briefly restate the model: the rate at which nest $m$ of site $x$ flips to state $i$, $c_i(x, m, \xi)$, is

$$c_0(x, m, \xi) = \begin{cases} 1, & \text{if } \xi^m(x) = 1, \\ 0, & \text{otherwise}, \end{cases}$$

$$c_1(x, m, \xi) = \begin{cases} \lambda n_1(x, \xi) n_2(x, \xi), & \text{if } \xi^m(x) = 0, \\ 0, & \text{otherwise}, \end{cases}$$

where $n_{m'}(x, \xi) = |\{z \in \mathcal{N}_x : \xi^{m'}(x+z) = 1\}|$ and $\mathcal{N}_x$ contains the site $x$ and its $2d$ nearest neighbors. The goal is to establish the existence of a phase transition.

2.1. *Existence of stationary distributions.* We first establish that stationary distributions exist. This is a generalization of Theorem 2.7 in Durrett [2] or Theorem III.2.3 in Liggett [5]. We use a method along the lines of Theorem 2.7 in Durrett [2]. Define

$$\bar{\xi}_0(x) = (1,1) \quad \text{for all } x$$

and let $\bar{\xi}_t$ be the process started with initial condition $\bar{\xi}_0$. Let $T_t f(\xi_0) = E^{\xi_0} f(\xi_t)$ be the semigroup corresponding to the particle system. $T_t$ is then a Feller semigroup by Theorem 1.1. We begin with a lemma.

LEMMA 2.1. *For any $A, B \subset S = \mathbb{Z}^d$, the function*

(13) $$t \mapsto P(\bar{\xi}_t^1(x) = 0 \ \forall x \in A, \bar{\xi}_t^2(y) = 0 \ \forall y \in B)$$

*is increasing.*

PROOF. Let $\alpha_0 = \bar{\xi}_s^1$ and $\beta_0 = \bar{\xi}_s^2$ for an arbitrary fixed $s$. Then $\bar{\xi}_0^1(x) \geq \alpha_0(x)$ and $\bar{\xi}_0^2(x) \geq \beta_0(x)$. Let $(\alpha_t, \beta_t)$ be the state at time $t$ of the particle system that started with initial condition $(\alpha_0, \beta_0)$. Then, by the fact that the particle system is monotone in initial conditions, we have

$$\bar{\xi}_t^1(x) \geq \alpha_t(x) \quad \text{and} \quad \bar{\xi}_t^2(x) \geq \beta_t(x)$$

for all $t$ and $x$. The Markov property of $\xi$ then implies that the function in (13) is increasing in $t$. □

THEOREM 2.2. *As $t \to \infty$, $\bar{\xi}_t$ converges weakly to $\bar{\xi}_\infty$. The limit is a stationary distribution that stochastically dominates all other stationary distributions and is called the* upper invariant measure.

PROOF. Let $A$ and $B$ be arbitrary subsets of $S$. For $C = \{x_1, \ldots, x_m\} \subset S$ and $D = \{y_1, \ldots, y_n\} \subset S$, we write

$$P(\bar{\xi}_t^1(z) = 0 \ \forall z \in A, \bar{\xi}_t^2(w) = 0 \ \forall w \in B, \bar{\xi}_t^1(x) = 1 \ \forall x \in C, \bar{\xi}_t^2(y) = 1 \ \forall y \in D)$$
$$= P(\bar{\xi}_t^1(z) = 0 \ \forall z \in A, \bar{\xi}_t^2(w) = 0 \ \forall w \in B) - P\left(\bigcup_{i=1}^{m+n} E_i\right),$$

where

$$E_i = \{\bar{\xi}_t^1(z) = 0 \ \forall z \in A \cup \{x_i\}, \bar{\xi}_t^2(w) = 0 \ \forall w \in B\} \quad \text{if } i = 1, \ldots, m$$

and

$$E_i = \{\bar{\xi}_t^1(z) = 0 \ \forall z \in A, \bar{\xi}_t^2(w) = 0 \ \forall w \in B \cup \{y_{i-m}\}\}$$
$$\text{if } i = m+1, \ldots, m+n.$$

We can use the inclusion-exclusion formula on $P(\bigcup_{i=1}^{m+n} E_i)$, that is,

$$P\left(\bigcup_{i=1}^{m+n} E_i\right) = \sum_{i=1}^{m+n} P(E_i) - \sum_{i<j} P(E_i \cap E_j)$$
$$+ \cdots + (-1)^{m+n+1} P(E_1 \cap \cdots \cap E_{m+n}).$$

Every term in the above expansion is in the form

$$P(\bar{\xi}_t^1(z) = 0 \ \forall z \in \cdot, \bar{\xi}_t^2(w) = 0 \ \forall w \in \cdot\cdot),$$

which is increasing in $t$ by Lemma 2.1. Therefore,

$$P(\bar{\xi}_t^1(z) = 0 \ \forall z \in A, \bar{\xi}_t^2(w) = 0 \ \forall w \in B, \bar{\xi}_t^1(x) = 1 \ \forall x \in C, \bar{\xi}_t^2(y) = 1 \ \forall y \in D)$$

converges for all $A$, $B$, $C$ and $D$ where $C$ and $D$ are finite, hence all finite-dimensional distributions converge. Thus, a weak limit $(\bar{\xi}_\infty^1, \bar{\xi}_\infty^2)$ exists and since $T_t$ is a Feller semigroup, $(\bar{\xi}_\infty^1, \bar{\xi}_\infty^2)$ is a stationary distribution. We can also easily see, via a monotonicity argument, that $(\bar{\xi}_\infty^1, \bar{\xi}_\infty^2)$ dominates all other stationary distributions. See the proof of Theorem 2.7 in Durrett [2] for a simpler version of this type of argument. □

2.2. *Extinction for sufficiently small $\lambda$.*

THEOREM 2.3. *If $\lambda|\mathcal{N}|^2 < 1$, then the particle system $\xi$ with generator (4) has no nontrivial stationary distribution.*

PROOF. We compare a modification of the particle system $\xi$ with the contact process. Theorem 2.6 of Durrett [2] states that if $\alpha|\mathcal{N}| < 1$, where $\alpha$ is the rate of infection, then the contact process has no nontrivial stationary distribution.

We modify the birth rates $c_1(x, m, \xi)$ in (2) and (3) of the particle model to obtain

$$c_1'(x, m, \xi) = \begin{cases} \lambda|\mathcal{N}|n_1(x, \xi), & \text{if } \xi^1(x) = 0, \\ \lambda|\mathcal{N}|n_2(x, \xi), & \text{if } \xi^2(x) = 0, \\ 0, & \text{otherwise}, \end{cases}$$

that is, birth of male (female) offspring no longer requires the presence of female (male) parents in the neighborhood. We denote this modified process $\tilde{\xi} = (\tilde{\xi}^1, \tilde{\xi}^2)$. The result of the modification is that $\tilde{\xi}^1$ and $\tilde{\xi}^2$ are now decoupled and $\tilde{\xi}^i$ behaves in the same way as the contact process with birth rate $\alpha = \lambda|\mathcal{N}|$. Furthermore, the modified process $(\tilde{\xi}^1, \tilde{\xi}^2)$ can also be constructed using the Poisson processes $\{T_n^{x,i,m} : n \geq 1\}$ and $\{U_n^{x,i,m} : n \geq 1\}$, just as in Section 1.1. Since $c_1'(x, m, \xi) \geq c_1(x, m, \xi)$ for all $x$, $m$ and $\xi$ [with the definition in either (2) or (3)], $(\tilde{\xi}^1, \tilde{\xi}^2)$ stochastically dominates the original process $(\xi^1, \xi^2)$. If $\alpha|\mathcal{N}| < 1$, then $\tilde{\xi}$ has no nontrivial stationary distribution and $(\tilde{\xi}_t^1, \tilde{\xi}_t^2)$ converges weakly to the all-0 state as $t \to \infty$ for any initial condition. Thus, $(\xi_t^1, \xi_t^2)$ also converges to the all-0 state for any initial condition if $\alpha|\mathcal{N}| = \lambda|\mathcal{N}|^2 < 1$, as required. □

2.3. *Survival for sufficiently large $\lambda$.* We use the idea of Chapter 4 of Durrett [2], that is, we compare the particle system to an oriented percolation process, see Durrett [2] for more details on the oriented percolation process. A particularly useful result that will be used in the proof below is Theorem 4.2 of Durrett [2].

THEOREM 2.4. *Let $W_n^p$ be an $M$-dependent oriented percolation process with density at least $1 - \gamma$, starting from the initial configuration $W_0^p$ in which the events $\{x \in W_0^p\}$, $x \in 2\mathbb{Z}$, are independent and have probability $p$. If $p > 0$ and $\gamma \leq 6^{-4(2M+1)^2}$, then*

$$\liminf_{n \to \infty} P(0 \in W_{2n}^p) \geq \tfrac{19}{20}.$$

This theorem shows that if the density of open sites $1 - \gamma$ is sufficiently close to 1 and we start with a Bernoulli initial condition for $W_0$, then the probability that 0 is wet at time $t$ does not go to 0 as $t \to \infty$, hence the upper invariant measure is nontrivial.

THEOREM 2.5. *If $\lambda$ is sufficiently large, then the particle system $\xi$ with generator (4) has a nontrivial stationary distribution.*

PROOF. We follow the method of proof as in Chapter 4 of Durrett [2]. First, we describe a construction of the particle system $\xi$ that is more specialized than the one given in Section 1.1. Recall that $S = \mathbb{Z}^d$ and $\mathcal{N} = \{x \in S : |x_1| + \cdots + |x_d| = 0 \text{ or } 1\}$. Let $m \in \{1, 2\}$, $x, y, z \in S$, $\{R_n^{x,m}, n \geq 1\}$ be independent Poisson processes with rate 1 and $\{T_n^{x,m,y,z}, n \geq 1\}$, with $y, z \in \mathcal{N}_x$, be independent Poisson processes with rate $\lambda$. At time $R_n^{x,m}$, any particle residing at $(x, m)$ is killed. And at time $T_n^{x,m,y,z}$, a particle is born at $(x, m)$ if $(x, m)$ is not already occupied and nests $(y, 1)$ and $(z, 2)$ are both occupied.

We will select an event $G_{\xi_0}$, measurable with respect to the filtration generated by all of the Poisson arrivals $\{R_n^{x,m}\}$ and $\{T_n^{x,m,y,z}\}$ at

$$x \in \{(-1, 0, \ldots, 0), (0, 0, \ldots, 0), (1, 0, \ldots, 0)\}$$

during the time interval $[0, t]$. For any $\gamma > 0$, regardless of how small, there exists $\lambda$ and $T$ and an event $G_{\xi_0}$ with

$$P(G_{\xi_0}) > 1 - \gamma,$$

so that on $G_{\xi_0}$, if $\xi_0(0,0,\ldots,0) = (1,1)$, then $\xi_T(1,0,\ldots,0) = \xi_T(-1, 0,\ldots,0) = (1,1)$. One can achieve this by choosing $T$ so small that the probability of any death occurring at any nests at sites $(-1,0,\ldots,0)$, $(0,0,\ldots,0)$ and $(1,0,\ldots,0)$ is less than $\gamma/2$; one can then choose $\lambda$ sufficiently large so that the probability of having birth events from $((0,0,\ldots,0),1)$ and $((0,0,\ldots,0),2)$ to each of the four nests at sites $(-1,0,\ldots,0)$ and $(1,0,\ldots,0)$ during $[0,T)$ is larger than $1 - \gamma/2$. In other words, if we define the event

$$\begin{aligned}G_{\xi_0} = \{&\text{There are no death events during } [0,T) \\ &\text{at sites } (-1,0,\ldots,0),\ (0,0,\ldots,0) \\ &\text{or } (1,0,\ldots,0);\ \text{and there are birth events from } ((0,0,\ldots,0),1) \\ &\text{and } ((0,0,\ldots,0),2) \text{ to each of the four nests at sites } (-1,0,\ldots,0) \\ &\text{and } (1,0,\ldots,0) \text{ during } [0,T)\},\end{aligned}$$

then $G_{\xi_0}$ satisfies the requirement and $P(G_{\xi_0}) > 1 - \gamma$ for some $\lambda$ and $T$. $G_{\xi_0}$ is the "good event" that will ensure male and female particles get born at sites $x - 1$ and $x + 1$ provided site $x$ is inhabited by both a male and a female particle.

Using this "good event" $G_{\xi_0}$, we can construct an oriented percolation process of density at least $1 - \gamma$ that is stochastically dominated by the particle system $\xi$, such that existence of nontrivial stationary distribution(s) for the oriented percolation process, as provided by Theorem 2.4, implies existence of nontrivial stationary distribution for the particle system. This part of the proof follows the proof of Theorem 4.3 in Durrett [2]. We therefore omit the details and instead refer interested readers to the proof of Theorem 3.3.2 in Yu [9]. □

**3. Convergence theorem for individual stirring.** In this section, we present the convergence result for the individual stirring model, as promised in Section 1.4. We work in a slightly more general setting and consider random processes

$$\xi_t^\varepsilon : \varepsilon \mathbb{Z}^d \times \{1, 2, \ldots, M\} \to \{0, 1, \ldots, \kappa - 1\}.$$

We call each $x \in \varepsilon \mathbb{Z}^d$ a site and each $(x, m) \in \varepsilon \mathbb{Z}^d \times \{1, 2, \ldots, M\}$ a nest. There are $M$ nests at each site. We think of the set of spatial locations $\mathbb{Z}^d \times \{1, 2, \ldots, M\}$ as consisting of $M$ floors of $\mathbb{Z}^d$. Let

$$\mathcal{N} = \{0, \varepsilon y_1, \ldots, \varepsilon y_N\}$$

be the interaction neighborhood of site 0. The process $\xi_t^\varepsilon$ evolves as follows.

1. *Birth and death.* The state of nest $(x, m)$ flips to $i$, $i = 0, \ldots, \kappa - 1$, at rate

$$c_i(x, m, \xi) = h_{i,m}(\xi(x, m), \xi(x + \varepsilon z_1, m_1), \ldots, \xi(x + \varepsilon z_L, m_L)),$$

where $L$ is a positive integer, $z_1, \ldots, z_L \in \mathcal{N}$, $m_1, \ldots, m_L \in \{1, 2, \ldots, M\}$ and

$$h_{i,m} : \{0, 1, \ldots, \kappa - 1\}^{L+1} \to \mathbb{R}^+$$

with $h_{i,m}(i, \ldots) = 0$.

2. *Rapid stirring.* For each $m \in \{1, 2, \ldots, M\}$ and $x, y \in \varepsilon \mathbb{Z}^d$ with $\|x - y\|_1 = \varepsilon$, $\xi^\varepsilon(x, m)$ and $\xi^\varepsilon(y, m)$ are exchanged at rate $\varepsilon^{-2}$.

This individual stirring model differs from the lily-pad stirring model described in Section 1.2 in that the stirring actions between corresponding nests at neighboring sites are now independent. More specifically, exchanges are allowed between neighboring nests on the same floor only, that is, between $(0, 1)$ and $(\varepsilon, 1)$, but not between $(0, 1)$ and $(\varepsilon, 2)$.

As an example, for $d = 1$, in the particle model with individual stirring with generator (7), we have $\kappa = 2$, $M = 2$, $L = 4$, $\mathcal{N} = \{0, -\varepsilon, \varepsilon\}$,

$$c_0(x, m, \xi) = \begin{cases} 1, & \text{if } \xi(x, m) = 1 \\ 0, & \text{otherwise,} \end{cases}$$

$$c_1(x, m, \xi) = \begin{cases} \lambda(\xi(x - \varepsilon, 1)\xi(x - \varepsilon, 2) + \xi(x + \varepsilon, 2)\xi(x + \varepsilon, 1)), \\ \quad \text{if } \xi(x, m) = 0, \\ 0, \quad \text{otherwise.} \end{cases}$$

In particular, we should define

$$(z_1, m_1) = (-1, 1), \qquad (z_2, m_2) = (-1, 2),$$
$$(z_3, m_3) = (1, 1), \qquad (z_4, m_4) = (1, 2)$$

and $h_i = h_{i,m}$ as

$$h_0(\alpha_0, \alpha_1, \alpha_2, \alpha_3, \alpha_4) = \alpha_0,$$
$$h_1(\alpha_0, \alpha_1, \alpha_2, \alpha_3, \alpha_4) = \lambda(\alpha_1 \alpha_2 + \alpha_3 \alpha_4)(1 - \alpha_0).$$

THEOREM 3.1. *Suppose $\{\xi_0^\varepsilon(x, m), (x, m) \in \varepsilon \mathbb{Z}^d \times \{1, 2, \ldots, M\}\}$ are independent and let $u_{i,m}^\varepsilon(t, x) = P(\xi_t^\varepsilon(x, m) = i)$. If $u_{i,m}^\varepsilon(0, x) = g_{i,m}(x)$ and $g_{i,m} : \mathbb{R}^d \to [0, 1]$ is $C^1$ for all $i$ and $m$ with $\sum_i g_{i,m} = 1$, then, for any smooth function $\phi$ with compact support, as $\varepsilon \to 0$,*

$$(14) \qquad \varepsilon^d \sum_{y \in \varepsilon \mathbb{Z}^d} \phi(y) 1_{\{\xi_t^\varepsilon(y,m) = i\}} \to \int \phi(y) u_{i,m}(t, y) \, dy \qquad \text{in probability,}$$

*where $u_{i,m}(t, x)$ is the bounded solution of*

$$\frac{\partial u_{i,m}}{\partial t} = \Delta u_{i,m} + f_{i,m}(u), \qquad u_{i,m}(0, x) = g_{i,m}(x),$$

$$f_{i,m}(u) = \langle c_i(0, m, \xi) 1(\xi(0, m) \neq i) \rangle_u - \sum_{j \neq i} \langle c_j(0, m, \xi) 1(\xi(0, m) = i) \rangle_u$$

and $\langle\phi(\xi)\rangle_u$ *denotes the expected value of* $\phi(\xi)$ *under the product measure in which state* $j$ *at nest* $m$ *has density* $u_{j,m}$, *that is,* $\xi(x,m)$, *with* $x \in \varepsilon\mathbb{Z}^d$ *and* $1 \leq m \leq M$, *are independent, with* $P(\xi(x,m)=j) = u_{j,m}$.

We refer interested readers to Theorem 4.0.5 of Yu [9] for its proof, which follows the approach used in the proof of Theorem 8.1 in Durrett [2]: first, a dual process is defined for the particle system, then this dual process is shown to be almost a branching random walk that converges to a branching Brownian motion as $\varepsilon \to 0$; furthermore, two different duals are asymptotically independent of each other. This asymptotic independence implies (14), as well as the convergence of the particle system itself to $u_{i,m}$.

In the proof of Theorem 3.1, the definition of the dual is the only part that differs from the proof of Theorem 8.1 in Durrett [2] in any significant way: each nest at a site requires a dual process and it is affected by any birth or death events that happen to any nests at that site, but dual processes for any two nests (even at the same site) are asymptotically independent.

**4. Invariant stationary distribution for lily-pad stirring.** In this section, we establish the existence of nontrivial stationary distribution(s) of the particle system with lily-pad stirring (as promised by Theorem 1.4) by showing that for sufficiently large $c$, the solution to

$$
\begin{aligned}
\frac{\partial u}{\partial t} &= \Delta u + (2c(1-u)+1)v - u, \\
\frac{\partial v}{\partial t} &= \Delta v + (c(u-v)-2)v,
\end{aligned}
\tag{15}
$$

with initial condition $u_0 = f, v_0 = g$, $f \geq g$ satisfies the following condition:

CONDITION (∗). There are constants $0 < D_1 < d_1 < d_2 < D_2 < 1$, $M$ and $T$ such that if $v_0(x) \in (D_1, D_2)$ for $x \in [-M, M]$, then $v_T(x) \in (d_1, d_2)$ for $x \in [-3M, 3M]$.

According to Chapter 9 of Durrett [2], this is a sufficient condition for the existence of nontrivial invariant stationary distribution(s) for the particle system with sufficiently fast stirring, so Theorem 1.4 will follow once Condition (∗) is established. Recall that Theorem 2.5 establishes that the dioecious particle model without rapid stirring has a nontrivial stationary distribution if the birth rate $\lambda$ is sufficiently large. If one works through the proof, however, one will find that "sufficiently large" in that argument means that $\lambda$ is larger than a number on the order of $6^{100}$, which is not too informative as to where exactly the critical $\lambda$ for the phase transition is. On the other hand, in the model with rapid stirring, one can obtain a far better idea of the range of $\lambda$ for which Condition (∗) holds.

In this work, we also establish Condition (∗) for sufficiently large $c$ (recall that $c = \lambda d$), but here, "sufficiently large" means that $c$ is only larger than a number on the order of $10^4$. We assume dimension $d = 1$; extension to $d > 1$ is straightforward. The proof consists of two parts: the first part, Section 4.1, establishes the existence of constants $d_1$ and $D_1$ and the second part, Section 4.2, establishes the existence of constants $d_2$ and $D_2$; the second part will be easy once the first part has been established.

4.1. *Lower bounds*: *Existence of $d_1$ and $D_1$ in Condition* (∗). Let

$$\mathcal{R} = \{(u,v) : 0 \leq v \leq u \leq 1\}$$

and $\eta$ be a vector field on $\mathcal{R}$ that generates a flow $(\mathcal{F}_\eta^t)_{t \geq 0}$ on $\mathcal{R}$, that is, if $(u_0, v_0) \in \mathcal{R}$, then $\mathcal{F}_\eta^t(u_0, v_0) \in \mathcal{R}$ for all $t \geq 0$. We assume $\mathcal{F}_\eta^t$ is monotone (in initial conditions), that is, it preserves the partial order on $\mathcal{R}$ given by

$$(u_1, v_1) \leq (u_2, v_2) \iff (u_1 \leq u_2 \text{ and } v_1 \leq v_2).$$

In other words, if $(u_1, v_1) \leq (u_2, v_2)$, then $\mathcal{F}_\eta^t(u_1, v_1) \leq \mathcal{F}_\eta^t(u_2, v_2)$. The scenario for $\eta$ that we consider is the following: the ODE system

$$\frac{du}{dt} = \eta_1(u, v), \qquad \frac{dv}{dt} = \eta_2(u, v)$$

has a stable fixed point $P_+$ close to the top corner of $\mathcal{R}$ with a relatively large basin of attraction, but $(0, 0)$ is also a stable fixed point (with a much smaller basin of attraction), and there is another unstable fixed point $P_-$ lying "between" these two stable ones; see Figure 1 for two examples of the vector field $\eta$ found in (15).

Theorem 9.2 in Durrett [2] establishes Condition (∗) for a specific predator-prey system with phase space $\{0, 1, 2\}$ at each site. The critical fact used in the proof is that the associated ODE system has only one interior equilibrium point and has a *global* Lyapunov function. The phase portrait of the ODE associated with (15), however, shows that it has two interior equilibrium points, $P_+$ and $P_-$, where $P_-$ is always a saddle point. Hence, we have a scenario where there is little hope of finding a global Lyapunov function and, in general, even finding an explicit Lyapunov function that works inside the domain of attraction of $P_+$ is difficult.

For general scalar reaction-diffusion equations (i.e., where the reaction terms are 1-dimensional), Chapter 15.4 of Taylor [7] provides an overview of methods and results. Convergence results in the scalar case, such as Condition (∗) in Chapter 9 of Durrett [2], can be established using techniques found in Fife and McLeod [4], just as in Durrett and Neuhauser [3]. For multidimensional reaction-diffusion equations, results regarding the existence of traveling wave solutions are more limited. Theorems 1.1, 1.2 and 4.2 of Chapter 3 of Volpert, Volpert and Volpert [8] provide existence results for traveling wave solutions for certain classes of monotone reaction-diffusion systems. Indeed, Theorem 1.1 of Chapter 3 of Volpert, Volpert and Volpert [8] applies to equation (13), but estimating the speed of the wave [which is essential for ensuring that Condition (∗) expands rather than shrinks] is still nontrivial and must be done on a case-by-case basis. An alternative approach to establishing Condition (∗) for equation (13) may be to use Theorem 1.1 of Volpert, Volpert and Volpert [8], which implies existence of traveling waves for (13), and then to try to obtain estimates for the speed of the wave (probably also a result involving sufficiently large $c$). With this estimate on the wave speed, one may then be able to use generalization of techniques in Fife and McLeod [4] to establish a convergence result.

The method we use to establish Condition (∗) for solutions of (15) is much more elementary and only requires the monotonicity property. Thus, it even applies to cases where the existence of traveling wave solutions is not known.

Consider the PDE system in one spatial dimension,

$$(16) \qquad \frac{\partial u}{\partial t} = \Delta u + \eta_1(u,v), \qquad \frac{\partial v}{\partial t} = \Delta v + \eta_2(u,v).$$

We first define the *shape* of the initial conditions $(u_0, v_0)$, which is a smoothed indicator function of a suitable interval. Let

$$(17) \qquad f_0(x) = \begin{cases} 1, & \text{if } x \in [-L+l, L-l], \\ h(x+L), & \text{if } x \in [-L-l, -L+l], \\ h(L-x), & \text{if } x \in [L-l, L+l], \\ 0, & \text{if } x \in (\infty, -L-l] \cup [L+l, \infty), \end{cases}$$

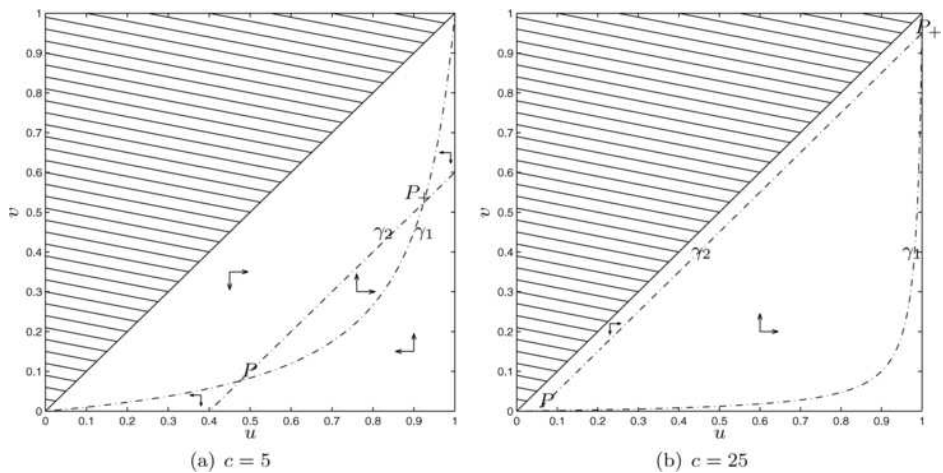

Fig. 1. *Phase space of the ODE:* $\eta_1 = 0$ *on* $\gamma_1$ *and* $\eta_2 = 0$ *on* $\gamma_2$.

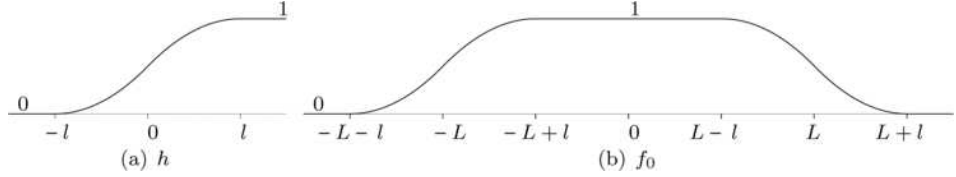

Fig. 2. *The functions $h$ and $f_0$.*

where $f_0 \in H^2(\mathbb{R})$ and

(18) $$h(x) = \begin{cases} 0, & x < -l, \\ \dfrac{1}{2}\left(\dfrac{x+l}{l}\right)^2, & -l \leq x \leq 0, \\ 1 - \dfrac{1}{2}\left(\dfrac{l-x}{l}\right)^2, & 0 < x \leq l, \\ 1, & x > l, \end{cases}$$

see Figure 2. In the above definition, the choice of $L$ is arbitrary, provided $L > l$, but later, we will choose $l$ small such that $|\Delta f_0|$ is large in $[-L-l, -L+l] \cup [L-l, L+l]$. We call the intervals $[-L-l, -L+l]$ and $[L-l, L+l]$ the *transition regions*. We observe that both $h$ and $h'$ are continuous at $x = 0$, with

$$h''(x) = \begin{cases} \dfrac{1}{l^2}, & \text{if } -l < x < 0, \\ -\dfrac{1}{l^2}, & \text{if } 0 < x < l, \end{cases}$$

so the graph of $h$ in the plane is symmetric about the point $(0, 1/2)$ and also

(19) $$|\Delta f_0| \leq \frac{1}{l^2}$$

everywhere.

We assume the solution $(u(t), v(t))$ to (16) starts with initial condition $(a_0 f_0, b_0 f_0)$. We would like to show that the interval over which the $v$-coordinate of the solution to (16) is $\geq b_0$ expands with time. Thus, intuitively, we would like to see the $v$-coordinate of $(u(t), v(t))$ increase, at least when $v(t)$ is larger than some threshold but smaller than $\pi_v(P_+)$, where $\pi_v(u', v') = v'$. If this were the case, we could use the $v$-coordinate as a "partial Lyapunov function" inside a subset of the basin of attraction of $P_+$. Unfortunately, this does not always hold, as is the case considered in this work where $\eta = ((2c(1-u)+1)v - u, (c(u-v)-2)v)$ and $\eta_2(u,v) < 0$ when $u - v$ is close to 0.

To overcome this difficulty, we define a convex family of nested subsets of $\mathcal{R}$ into which the flow $\mathcal{F}_\eta$ contracts at a sufficiently large rate. The upper boundary of these nested subsets consists of lines $u = v$ and $u = 1$, while the lower boundary consists of curves from the family $\{\gamma_\theta\}_{\theta \in [0, \theta_0]}$, where $\gamma_\theta = (1+\theta)\gamma \cap \mathcal{R}$ and $\gamma$ is a nonincreasing (i.e., with nonincreasing $v$-coordinate when parameterized according to increasing $u$-coordinate) curve whose exact definition may change depending on the exact form of $\eta$. Assume, furthermore, that there is a family of continuous mappings $\varphi_\theta : [0,\theta_0] \times \gamma_\theta \to \gamma_\theta$ such that $\mathcal{F}_\eta$ satisfies the following.

ASSUMPTION 4.1.  There exist $\varepsilon, K_1, K_2 > 0$ ($K_1$ is large and $K_2$ is small) and sufficiently small $s_0$ such that for $\theta \in [0,\theta_0]$, $(a_0,b_0) \in \gamma_\theta$ and $s \in [0,s_0]$, we have

$$\mathcal{F}_\eta^s(\alpha a_0, \alpha b_0) \geq \begin{cases} ((1+K_1 s)\alpha a_s, (1+K_1 s)\alpha b_s), & \text{if } \alpha b_0 \geq \varepsilon, \\ ((1-K_2 s)\alpha a_s, (1-K_2 s)\alpha b_s), & \text{if } \alpha b_0 < \varepsilon, \end{cases}$$

where $(a_s, b_s) = \varphi_\theta(s, (a_0, b_0))$.

For each $s$, the mapping $\varphi_\theta(s, \cdot)$ maps points on the curve $\gamma_\theta$ to other points still on $\gamma_\theta$. Under $\mathcal{F}_\eta$, the upper part (the part above the horizontal line $v = \varepsilon$) of the line segment $\{(\alpha a_0, \alpha b_0) : 0 \leq \alpha \leq 1\}$ is pulled above $(\alpha a_s, \alpha b_s)$, while its lower part is not pushed too far beneath $(\alpha a_s, \alpha b_s)$; see Figure 3 for a schematic illustration. Note that $((1+K_1 s)a_s, (1+K_1 s)b_s)$ again lies on a curve $\gamma_{\bar{\theta}}$, where $\bar{\theta} > \theta + \delta s$ and $\delta$ depends on $K_1$, $K_2$ and the geometry of $\gamma$.

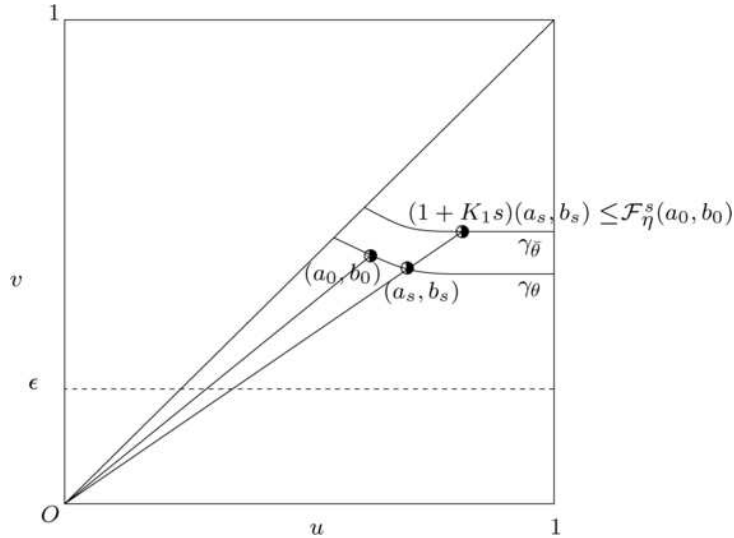

FIG. 3.  *Illustration of Assumption 4.1.*

PROPOSITION 4.2. *Assume that $\mathcal{F}_\eta$ satisfies Assumption* 4.1 *with* $\min\{v:(u,v) \in \gamma\} \geq 0.5$, $\varepsilon \leq 0.24$ *and*

$$\text{(20)} \qquad \frac{K_1}{K_2} > 21\left(\frac{200}{199}\right)^2.$$

*There then exists $l > 0$ [as in the definition of $h$ in (18)], $\delta_1$, $\delta_2$, $s_0 > 0$ such that for $\theta \in [0, \theta_0]$, $(a_0, b_0) \in \gamma_\theta$ and $s \in [0, s_0]$,*

$$e^{s\Delta}\mathcal{F}_\eta^s(a_0 f_0, b_0 f_0) \geq ((1+\delta_2 s)a_s f_s, (1+\delta_2 s)b_s f_s),$$

*where $(a_s, b_s) = \varphi_\theta(s, (a_0, b_0))$ as above, $f_0$ is defined as in (17) and $f_s$ is $f_0$ with the transition regions translated outward by $\delta_1 s$:*

$$\text{(21)} \quad f_s(x) = \begin{cases} 1, & x \in [-L+l-\delta_1 s, L-l+\delta_1 s], \\ h(x + L + \delta_1 s), & x \in [-L-l-\delta_1 s, -L+l-\delta_1 s], \\ h(L + \delta_1 s - x), & x \in [L-l+\delta_1 s, L+l+\delta_1 s], \\ 0, & x \in (\infty, -L-l-\delta_1 s] \cup [L+l+\delta_1 s, \infty). \end{cases}$$

REMARK 4.3. By abuse of notation, we again let $\mathcal{F}_\eta$ denote the time evolution of the spatial system where each $(u(x), v(x))$, $x \in \mathbb{R}$, flows independently along the vector field $\eta$.

For any $(u, v) \in \mathcal{R}$, the set of values $\{(uf_0(x), vf_0(x)): x \in \mathbb{R}\}$ forms a line segment in $\mathcal{R}$ with endpoints $O$ and $(u, v)$. Proposition 4.2 above combines the properties of the flow $\mathcal{F}_\eta^s$ for small $s$ with the spatial distribution generated by the heat flow. Before we prove Proposition 4.2, we state two technical lemmas necessary for its proof.

LEMMA 4.4. *If $l$ is fixed and $f = f_0$ is defined as in (17), then for*

$$x \in \left(-L-l-s, -L-\frac{l}{200}\right) \cup \left(L+\frac{l}{200}, L+l+s\right)$$

*and $s$ small, we have*

$$e^{s\Delta}f(x) \geq f(x) + \frac{s}{5l^2}.$$

LEMMA 4.5. *Let $s > 0$ be fixed, $f_0$ as defined in (17) and*

$$\hat{f}(x) = \begin{cases} f_0(x) + ms, & -L-l-s < x < L+l+s, \\ 0, & otherwise, \end{cases}$$

*where $m > 0$. Then there exist positive constants $\delta_2$ depending on $m$ but independent of $s$ such that for all $x$, $\hat{f}(x) \geq (1+\delta_2 s)f_s(x)$, where $f_s$ is defined as in (21).*

Lemma 4.4 says that the lower part of the transition region in $f_0$ increases at a rate proportional to $1/l^2$, which is the case since in that part of the transition region, $f_0$ is convex and $\Delta f_0 = O(1/l^2)$. It can be proved by working with the solution to the heat equation in convolution form $e^{s\Delta} g(x) = \frac{1}{\sqrt{4\pi s}} \int \exp(-y^2/4s) g(x-y) \, dy$, or working with Brownian motions if one prefers to calculate probabilities. Lemma 4.5 involves only elementary calculus. We refer interested readers to Yu [9] (Lemmas 5.1.6 and 5.1.7) for full details of proofs of these two lemmas.

PROOF OF PROPOSITION 4.2. By virtue of (20), we can choose $l > 0$ such that

$$K_1 > \left(\frac{200}{199}\right)^2 \frac{4.02}{l^2}, \qquad K_2 < \frac{1}{5.05 l^2}.$$

Hence, we can select $m > 0$ such that

$$(22) \qquad m \geq \frac{K_1}{2}\left(\frac{199}{200}\right)^2 - \frac{2.01}{l^2}$$

and

$$(23) \qquad m \geq \frac{1}{5l^2} - 1.01 K_2.$$

Fix $\theta \in [0, \theta_0]$, $(a_0, b_0) \in \gamma_\theta$ and $s \in [0, s_0]$ for the moment.

First, consider $x \in [-L - \frac{1}{200}, L + \frac{1}{200}]$, where the definition of $f_0$ in (17) implies that

$$a_0 f_0(x) \geq b_0 f_0(x) \geq \tfrac{1}{2} f_0(L + \tfrac{1}{200}) = \tfrac{1}{2}\tfrac{1}{2}(\tfrac{199}{200})^2 \geq 0.24 \geq \varepsilon.$$

Hence, by the first half of Assumption 4.1,

$$\mathcal{F}_\eta^s(a_0 f_0(x), b_0 f_0(x))$$
$$\geq ((1 + K_1 s) a_s f_0(x), (1 + K_1 s) b_s f_0(x))$$

for $s \in [0, s_0]$ and $x \in [-L - \frac{1}{200}, L + \frac{1}{200}]$. Furthermore, by (19), we have $e^{s\Delta} f_0(x) \geq f_0(x) - 2s/l^2$ for all $x$ and $s \geq 0$. This, together with the monotonicity of $e^{s\Delta}$, shows that

$$(e^{s\Delta}\mathcal{F}_\eta^s)(a_0 f_0, b_0 f_0)(x)$$
$$\geq e^{s\Delta}((1 + K_1 s) a_s f_0(x), (1 + K_1 s) b_s f_0(x))$$
$$\geq \left((1 + K_1 s) a_s\left(f_0(x) - \frac{2s}{l^2}\right), (1 + K_1 s) b_s\left(f_0(x) - \frac{2s}{l^2}\right)\right).$$

Using the fact that $f_0(x) \geq f(L + \frac{1}{200}) = \frac{1}{2}(\frac{199}{200})^2$ for $x \in [-L - \frac{1}{200}, L + \frac{1}{200}]$, we can bound the $v$-coordinate in the above as follows:

$$b_s(1 + K_1 s)\left(f_0(x) - \frac{2s}{l^2}\right) = b_s\left(f_0(x) + K_1 f_0(x) s - \frac{2s}{l^2} - \frac{2K_1 s^2}{l^2}\right)$$

$$\geq b_s\left[f_0(x) + \frac{K_1}{2}\left(\frac{199}{200}\right)^2 s - \frac{2s}{l^2} - \frac{2K_1 s^2}{l^2}\right]$$

$$\geq b_s\left[f_0(x) + \left(\frac{K_1}{2}\left(\frac{199}{200}\right)^2 - \frac{2.01}{l^2}\right)s\right],$$

if $s$ is sufficiently small. The $u$-coordinate can be treated analogously, therefore

$$(24) \quad (e^{s\Delta}\mathcal{F}^s_\eta)(a_0 f_0, b_0 f_0)(x) \geq (a_s(f_0(x) + ms), b_s(f_0(x) + ms)),$$

by (22).

We now consider $x \in [-L - l - s, -L - \frac{1}{200}] \cup [L + \frac{1}{200}, L + l + s]$. By the second half of Assumption 4.1, we have

$$\mathcal{F}^s_\eta(\alpha a_0, \alpha b_0) \geq ((1 - K_2 s)a_s f_0(x), (1 - K_2 s)b_s f_0(x))$$

for such $x$. Combining this with Lemma 4.4 and the monotonicity of $e^{s\Delta}$, we obtain

$$(e^{s\Delta}\mathcal{F}^s_\eta)(a_0 f_0, b_0 f_0)(x) \geq e^{s\Delta}((1 - K_2 s)a_s f_0(x), (1 - K_2 s)b_s f_0(x))$$

$$\geq \left((1 - K_2 s)a_s\left(f_0(x) + \frac{s}{5l^2}\right), (1 - K_2 s)b_s\left(f_0(x) + \frac{s}{5l^2}\right)\right)$$

$$(25)$$

$$\geq \left(a_s\left[f_0(x) + \left(\frac{1}{5l^2} - 1.01K_2\right)s\right], b_s\left[f_0(x) + \left(\frac{1}{5l^2} - 1.01K_2\right)s\right]\right)$$

$$\geq (a_s(f_0(x) + ms), b_s(f_0(x) + ms))$$

by (23), where the third inequality requires $s$ to be sufficiently small. Combining (24), (25) and Lemma 4.5 then yields the claim. □

COROLLARY 4.6. *Suppose that $\eta$ satisfies Assumption 4.1 with* $\min\{v : (u,v) \in \gamma\} = D_1 < d_1 < \min\{v : (u,v) \in (1+\theta_0)\gamma\}$. *If $T$ is sufficiently large and $v_0(x) > D_1$ for $x \in [-L+l, L-l]$, then $v_T(x) > d_1$ for $x \in [-3L, 3L]$.*

PROOF. This is an easy consequence of Proposition 4.2 and the nonlinear Trotter product formula (Proposition 15.5.2 from Taylor [7])

$$(u_t, v_t) = \lim_{n \to \infty} (e^{(t/n)\Delta} \mathcal{F}^{t/n})^n(f, g).$$

Here, the convergence occurs in the space $\mathcal{BC}^1(\mathbb{R})$, which consists of functions $f$ such that $f'$ is bounded and continuous on $\mathbb{R}$ and both $f$ and $f'$ can be extended to be continuous on the compactification $\widehat{\mathbb{R}}$ via the point at infinity; the norm used here is $\|\cdot\|_\infty + \|\partial/\partial x(\cdot)\|_\infty$. For sufficiently large $n$, Proposition 4.2 says that application of $e^{(1/n)\Delta}\mathcal{F}^{1/n}$ to the function $(a_0 f_0, b_0 f_0)$ with $(a,b) \in \gamma_\theta$ yields a result that is larger than $(\tilde{a}\tilde{f}, \tilde{b}\tilde{f})$, where $(\tilde{a}, \tilde{b}) \in \gamma_{\theta+\delta/n}$ for some $\delta > 0$ (whose exact value depends on $\delta_2$ and $\theta$) and the "flat region" in $\tilde{f}$ is at least $2\delta_1/n$ larger than that of $f_0$. $\square$

It remains to check that if $c$ is sufficiently large, then the vector field

$$
\begin{aligned}
\eta_1(u,v) &= (2c(1-u)+1)v - u, \\
\eta_2(u,v) &= (c(u-v)-2)v
\end{aligned}
\tag{26}
$$

as in (15) satisfies Assumption 4.1.

LEMMA 4.7. *If $c$ is sufficiently large, then we can find $s_0$, $\theta_0$, $K_1$, $K_2$, $\varepsilon > 0$ and a curve $\gamma$ such that Assumption* 4.1 *and the assumptions of Proposition* 4.2 *are satisfied for the flow generated by* (26).

PROOF. For $c \geq 8{,}800$, we construct a vector field $\xi = (\xi_1(u,v), \xi_2(u,v))$ for $(u,v) \in \mathcal{R}$ such that

$$\xi_1(u,v) \leq \eta_1(u,v), \qquad \xi_2(u,v) \leq \eta_2(u,v), \qquad (u,v) \in \mathcal{R}$$

and we will show that the flow $\mathcal{F}_\xi$ generated by $\xi$ satisfies Assumption 4.1 and the assumptions of Proposition 4.2 with $\varepsilon = 0.24$, $K_1 = 44$ and $K_2 = 2$. Note that there is a fair amount of leeway in the choice of the constants, as we have not striven for optimality in that respect.

Define $A = (0.51, 0.51)$, $B = (0.55, 0.5)$, $C = (0.9, 0.5)$, $D = (1.0, 0.5)$ and let the curve $\gamma$ be given be $\gamma = \gamma_0 = \overline{AB} \cup \overline{BD}$ (see Figure 4), where $\overline{AB}$ is the line segment connecting $A$ and $B$. For $\theta \in [-0.54, 0.2]$, let

$$
\begin{aligned}
A_\theta &= (1+\theta)A, & B_\theta &= (1+\theta)B, \\
C_\theta &= (0.9, (1+\theta)0.5), & D_\theta &= (1.0, (1+\theta)0.5)
\end{aligned}
$$

and

$$\gamma_\theta = \overline{A_\theta B_\theta} \cup \overline{B_\theta C_\theta}.$$

Hence, $\gamma_\theta = (1+\theta)\gamma \cap \mathcal{R}$. Put

$$
\begin{aligned}
A' &= A_{-0.54}, & B' &= B_{-0.54}, & C' &= C_{-0.54}, & D' &= D_{-0.54}, \\
A'' &= A_{0.2}, & B'' &= B_{0.2}, & C'' &= C_{0.2}, & D'' &= D_{0.2}.
\end{aligned}
$$

If we let $\varepsilon' = 0.23$, then $B'$, $C'$ and $D'$ lie on the horizontal line $\{v = \varepsilon'\}$, while $B''$, $C''$ and $D''$ lie on the horizontal line $\{v = 0.6\}$ and $\pi_v(A'') = 0.2346 < \varepsilon$.

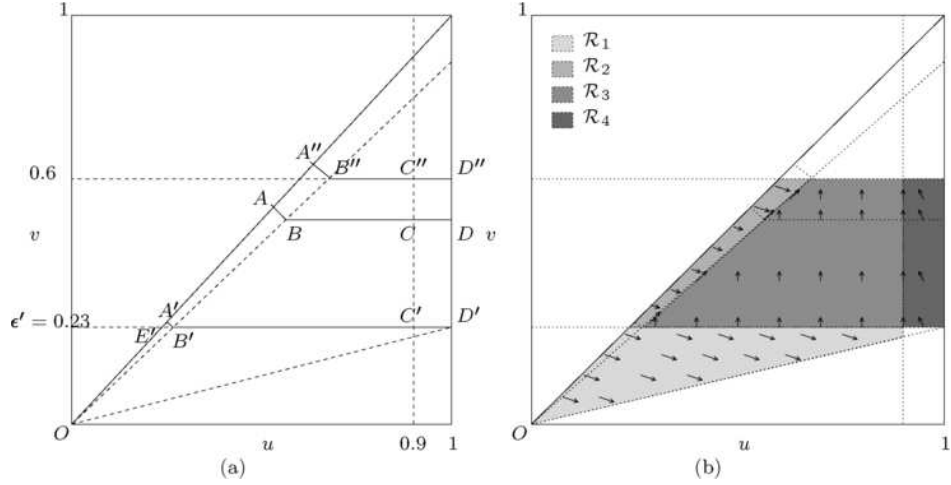

FIG. 4. (a) *sketch of the curves* $\gamma_\theta$. (b) *sketch of* $\xi$, *arrow sizes and directions not drawn to scale*.

We define the vector field $\xi$ as a piecewise linear function in the following way: let $F_1 = 400$ and $F_2 = 75$. This implies that

$$\frac{1}{2\sqrt{17}}((F_1 - 8) \wedge 5.1 F_2) \geq K_1, \tag{27}$$

which will be needed later on. For $(u, v) \in \mathcal{R}$, let $\xi = (\xi_1(u, v), \xi_2(u, v))$ be given by

$$\begin{aligned}
(F_1 u, -2v) \quad & \text{in } \mathcal{R}_1 = \{\varepsilon' u < v \leq \varepsilon', u \leq 0.9\} \\
& \text{or } \mathcal{R}_2 = \{\varepsilon' < v \leq 0.6, u < 1.1 v\}, \\
(1.1 F_2 v, F_2 v) \quad & \text{on } L_1 = \overline{B' B''}, \\
(0, F_2 v) \quad & \text{in } \mathcal{R}_3 = (\text{interior of the trapezoid } B' C' C'' B'') \cup \overline{C' C'''}, \\
(-u, F_2 v) \quad & \text{in } \mathcal{R}_4 = \text{interior of rectangle } C' D' D'' C''', \\
(\eta_1(u, v), \eta_2(u, v)) \quad & \text{in the rest of } \mathcal{R};
\end{aligned}$$

see Figure 4. To verify that $\xi \leq \eta$ on $\mathcal{R}$ if $c \geq 8,800$, we first observe that $\eta_1 \geq -u$ and $\eta_2 \geq -2v$ are trivial bounds that hold for all of $\mathcal{R}$. The other definitions of $\xi$ can be verified region by region, as follows.

$\mathcal{R}_1 \cup \mathcal{R}_2$: Since $v \geq \varepsilon' u$ and $u \leq 0.9$, we have

$$\eta_1 \geq (0.2c + 1)\varepsilon' u - u \geq (0.046 c - 1) u \geq 400 u,$$

if $c \geq 8,800$.

$L_1$: The validity of $\eta_1 \geq \xi_1$ comes from the same argument as for $\mathcal{R}_1 \cup \mathcal{R}_2$. For $\xi_2$, since $v \geq \varepsilon' = 0.23$ and $u = 1.1v$, we have

$$\eta_2(u,v) = (0.1cv - 2)v \geq 75v,$$

if $c \geq 3,400$.

$\mathcal{R}_3$: Since $u \leq 0.6$ and $v \geq \varepsilon' u$, we have

$$\eta_1(u,v) \geq (2c(0.4) + 1)\varepsilon' u - u,$$

if $c \geq 5$. The validity of $\eta_2 \geq \xi_2$ comes from the same argument as for $L_1$.

$\mathcal{R}_4$: The validity of $\eta_2 \geq \xi_2$ comes from the same argument as for $L_1$.

Let $s_0 = [(\log(12/11))/F_2] \wedge [(\log(45/23))/(F_1 + 2)]$. These choices guarantee that for $s \in [0, s_0]$, $\pi_v(\mathcal{F}_\xi^s(u_0, v_0)) \leq 0.6$ if $v_0 \leq 0.55$ and

$$\mathcal{F}_\xi^s(u_0, v_0) = (e^{F_1 s} u_0, e^{-2s} v_0) \in \{v \geq \tfrac{23}{90} u\}$$

if $v_0 \geq 0.5 u_0$. For $\theta \in [0, 0.2]$, $s \in [0, s_0]$ and $(u_0, v_0) \in \gamma_\theta$, define

(28)　$(u_s, v_s) = \varphi_\theta(s, (u_0, v_0)) = $ intersection of $\overline{O \mathcal{F}_\xi^s(u_0, v_0)}$ and $\gamma_\theta$.

Let $\Delta A'B'E'$ denote the triangle with vertices $A'$, $B'$ and $E' = (\varepsilon', \varepsilon')$ (a subset of $\mathcal{R}_2$). We claim that

(29)　$\alpha(u_0, v_0) \notin \mathcal{R}_1 \cup \Delta A'B'E' \implies \mathcal{F}_\xi^s(\alpha u_0, \alpha v_0) \geq \alpha(1 + K_1 s)(u_s, v_s),$

(30)　$\alpha(u_0, v_0) \in \mathcal{R}_1 \cup \Delta A'B'E' \implies \mathcal{F}_\xi^s(\alpha u_0, \alpha v_0) \geq \alpha(1 - K_2 s)(u_s, v_s).$

Notice that this implies the required inequality in Assumption 4.1 because $\mathcal{R}_1 \cup \Delta A'B'E'$ lies beneath the horizontal line $v = \varepsilon = 0.24$.

PROOF OF (29). First, consider the case $(u_0, v_0) \in \overline{A_\theta B_\theta}$. As $\xi$ is linear in $\mathcal{R}_1$ and $\mathcal{R}_2$, it suffices to restrict to $\theta = 0$ and $\alpha = 1$ in this case, that is, consider $(u_0, v_0) = \lambda A + (1 - \lambda)B \in \gamma$ for some $\lambda \in [0, 1]$. Note that $(1/\sqrt{17}, 4/\sqrt{17})$ is a unit normal vector perpendicular to $\overline{AB}$. Provided $\mathcal{F}_\xi^s(u_0, v_0)$ is to the left of $\overline{B'B''}$, the rate of increase of $\mathcal{F}_\xi^s(u_0, v_0) - (u_0, v_0)$ in the direction of $(1/\sqrt{17}, 4/\sqrt{17})$ is

$$(F_1 \pi_u(\mathcal{F}_\xi^s(u_0, v_0)), -2\pi_v(\mathcal{F}_\xi^s(u_0, v_0))) \cdot (1/\sqrt{17}, 4/\sqrt{17})^t$$
$$\geq \frac{(F_1 - 8)\pi_v(\mathcal{F}_\xi^s(u_0, v_0))}{\sqrt{17}} \geq \frac{(F_1 - 8)0.5}{\sqrt{17}} \geq K_1,$$

by (27), where we can use the lower bound $\pi_v(\mathcal{F}_\xi^s(u_0, v_0)) \geq 0.5$ since the trajectory $\mathcal{F}_\xi^s(u_0, v_0)$ for $(u_0, v_0) \in \overline{A_\theta B_\theta}$ stays above $\overline{A_\theta B_\theta}$. On the other hand, once the trajectory $\mathcal{F}_\xi^s(u_0, v_0)$ hits $\overline{B'B''}$, the rate of increase of $\mathcal{F}_\xi^s(u_0, v_0) - (u_0, v_0)$ in the direction of $(1/\sqrt{17}, 4/\sqrt{17})$ is

$$\pi_v(\mathcal{F}_\xi^s(u_0, v_0))(1.1F_2, F_2) \cdot (1/\sqrt{17}, 4/\sqrt{17})^t$$
$$= \frac{5.1 F_2 \pi_v(\mathcal{F}_\xi^s(u_0, v_0))}{\sqrt{17}} \geq \frac{5.1 F_2 \cdot 0.5}{\sqrt{17}} \geq K_1,$$

by (27). Therefore,

$$\mathcal{F}_\xi^s(u_0, v_0) \geq (1 + K_1 s)(u_s, v_s),$$

if $(u_0, v_0) \in \overline{A_\theta B_\theta}$.

Now, consider $(u_0, v_0) \in \overline{B_\theta C_\theta}$, with $\theta \in [0, 0.2]$ and $\alpha \in (0, 1]$ such that $\alpha v_0 \geq \varepsilon'$. Then $(u_0, v_0)$ and $(\alpha u_0, \alpha v_0)$ are in $\mathcal{R}_3$. Hence,

$$\mathcal{F}_\xi^s(\alpha u_0, \alpha v_0) = (\alpha u_0, \alpha e^{F_2 s} v_0)$$

until, possibly, the trajectories hit $\overline{B'B''}$ (which will happen simultaneously by the linearity of $\xi$ in $\mathcal{R}_3$), when the $u$-coordinate will also start to increase. Furthermore, for $(u_0, v_0) \in \overline{B_\theta D_\theta}$, $v_s = v_0$ by the geometry of $\gamma_\theta$ and the definition of $\varphi$ in (28), thus

$$\mathcal{F}_\xi^s(\alpha u_0, \alpha v_0) \geq \alpha(1 + K_1 s)(u_s, v_s)$$

in this case as well since $F_2 > K_1$. Similar reasoning applies when $(u_0, v_0) \in \overline{C_\theta D_\theta}$:

$$\mathcal{F}_\xi^s(\alpha u_0, \alpha v_0) \geq (\alpha u_0 e^{-s}, \alpha v_0 e^{F_2 s}),$$

where we have equality provided $\mathcal{F}_\xi^s(\alpha u_0, \alpha v_0) \in \mathcal{R}_4$. Thus, before $\mathcal{F}_\xi(u_0, v_0)$ hits $\overline{C'C''}$, we have $(u_s, v_s) = (u_0 e^{-(F_2+1)s}, v_0)$ and

$$\mathcal{F}_\xi^s(\alpha u_0, \alpha v_0) \geq \alpha e^{F_2 s}(u_s, v_s) \geq \alpha(1 + K_1 s)(u_s, v_s).$$

Notice that if $(\alpha u_0, \alpha v_0) \in \mathcal{R}_4$, then the $u$-coordinates of $\mathcal{F}_\xi^s(\alpha u_0, \alpha v_0)$ will stop decreasing once $\overline{C'C''}$ is hit. □

PROOF OF (2.5). According to the definition of $\xi$ in $\mathcal{R}_1 \cup \mathcal{R}_2$,

(31) $$\mathcal{F}_\xi^s(\alpha u_0, \alpha v_0) = \alpha(u_0 e^{F_1 s}, v_0 e^{-2s})$$

for $\alpha(u_0, v_0) \in \mathcal{R}_1 \cup \Delta A'B'E'$, provided the trajectories remain in $\mathcal{R}_1 \cup \Delta A'B'E'$, which is the case by the choice of $s_0$. Assume, first, that $(u_0, v_0) \in \overline{A_\theta B_\theta}$. If $\mathcal{F}_\xi(u_0, v_0)$ does not hit $\overline{B'B''}$ [i.e., $\mathcal{F}_\xi(\alpha u_0, \alpha v_0)$ does not hit $\overline{OB'}$] by time $s$, then

$$\mathcal{F}_\xi^s(\alpha u_0, \alpha v_0) = \alpha \mathcal{F}_\xi^s(u_0, v_0),$$

hence, in particular, even

$$\mathcal{F}^s_\xi(\alpha u_0, \alpha v_0) \geq \alpha(u_s, v_s),$$

which implies (30). On the other hand, if $\mathcal{F}_\xi(u_0, v_0)$ does hit $\overline{B'B''}$ by time $s$, then $(u_s, v_s) = \varphi_\theta(s, (u_0, v_0)) = B_\theta$ and, hence,

$$\mathcal{F}^s_\xi(\alpha u_0, \alpha v_0) \geq \alpha(u_s, v_s e^{-2s}) \geq \alpha(1 - K_2 s)(u_s, v_s),$$

which establishes (30).

Now, consider the case $(u_0, v_0) \in \overline{B_\theta D_\theta}$, where we always have $(u_s, v_s) \leq (u_0, v_0)$ ($u_s$ decreases with $s$ until the trajectory hits $\overline{B'B''}$, after which both $u_s$ and $v_s$ remain constant). Thus, (31) implies

$$\mathcal{F}^s_\xi(\alpha u_0, \alpha v_0) \geq \alpha e^{-2s}(u_0, v_0) \geq \alpha(1 - K_2 s)(u_s, v_s),$$

so (30) also holds in this case. The proof of Lemma 4.7 is now complete. $\square$

4.2. *Upper bounds*: *Existence of $d_2$ and $D_2$ in Condition (∗)*. We first establish the following proposition.

PROPOSITION 4.8. *If $c$ is sufficiently large, then there exist constants $d_2 < D_2 < 1$ and $T$ such that if $v_0(x) < D_2$ for $x \in [-L + l, L - l]$, then $v_t(x) < d_2$ for $x \in [-3L, 3L]$ for all $t \geq T$, where $(u_t, v_t)$ solves the PDE (15).*

PROOF. Because of the monotonicity of (15), it suffices to consider the uniform initial condition $u_0 \equiv 1, v_0 \equiv 1$ and to show that $v_t < d_2$ for large $t$. Therefore, we need only concern ourselves with the ODE

$$\frac{du}{dt} = (2c(1-u) + 1)v - u,$$

$$\frac{dv}{dt} = (c(u-v) - 2)v.$$

We can bound $\eta_2(u, v) = (c(u-v) - 2)v$ for any $v > 1 - 1/c$ as follows:

$$\eta_2(u, v) = (c(u-v) - 2)v \leq \left(c\left(1 - \left(1 - \frac{1}{c}\right)\right) - 2\right)v$$

$$= -v < -\left(1 - \frac{1}{c}\right) < 0,$$

if $c > 1$. Thus, for any $d_2$ satisfying $1 > d_2 > 1 - 1/c$, there exists $T$, such that if $u_0 = v_0 \equiv 1$, then $v_t < d_2$ for $t \geq T$. $\square$

PROOF OF THEOREM 1.4. If we let $M = L - l$, then Corollary 4.6, Lemma 4.7 and Proposition 4.8 show that Condition (∗) from the beginning of Section 4 holds for the PDE system (15). This, in turn, implies the conclusion of the theorem by section 3 of Durrett and Neuhauser [3]. $\square$

**Acknowledgments.** The author is especially grateful to Edwin Perkins for his generous help and support throughout this research and for reading the manuscript multiple times. The author is also grateful to Martin Barlow for helpful discussions and suggestions. In addition, the author wishes to thank the two anonymous referees for valuable suggestions and comments, one of whom wrote seven page improved version of the proof for establishing Condition (∗) for lily-pad stirring; Section 4 in this work is based on that referee's rewrite of the original proof, which was less streamlined and harder to read than the current version.

Department of Statistics
University of Oxford
1 South Parks Road, Oxford OX1 3TG
United Kingdom
E-mail: fyz@stats.ox.ac.uk
URL: http://www.stats.ox.ac.uk/~fyz